\documentclass[a4paper,12pt,oneside]{article}
\usepackage[english]{babel}
\usepackage[T1]{fontenc} 
\usepackage[utf8]{inputenc}
\usepackage{amsthm}
\usepackage{amsmath}
\usepackage{amssymb}  
\usepackage{indentfirst}
\usepackage{fancyhdr}
\usepackage{amsthm}
\usepackage{graphicx}
\usepackage{pdfpages}
\usepackage{esint}
\usepackage{url}

\definecolor{dg}{rgb}{0.01, 0.75, 0.24}

\numberwithin{equation}{section}

\theoremstyle{plain} 
\newtheorem{thm}{Theorem}[section]

\newtheorem{lem}[thm]{Lemma}

\def\io{\int_{\Omega}}
\def\R{\mathbb{R}}
\def\N{\mathbb{N}}
\def\K{\mathcal{K}_{\psi}(\Omega)}
\def\dd{\mathrm{d}}
\def\supp{\mathrm{supp}}

\def\fibr2 {\displaystyle\fint_{B_{\frac{R}{2}}}}
\def\fibro2 {\fint_{B_{\frac{\rho}{2}}}}

\usepackage{geometry}
\allowdisplaybreaks[4]

\begin{document}

\title{A local boundedness result for a class of obstacle problems with non-standard growth conditions}

\author {   
\sc{Mariapia De Rosa}\thanks{Dipartimento di Matematica e Applicazioni "R. Caccioppoli", Università degli Studi di Napoli "Federico II", Via Cintia, 80126 Napoli
 (Italy). E-mail: \textit{mariapia.derosa@unina.it}} 
  \;\textrm{and}   Antonio Giuseppe Grimaldi \thanks{Dipartimento di Matematica e Applicazioni "R. Caccioppoli", Università degli Studi di Napoli "Federico II", Via Cintia, 80126 Napoli (Italy).
E-mail: \textit{antoniogiuseppe.grimaldi@unina.it}}}

\maketitle
\maketitle

\begin{abstract}We prove the local boundedness for solutions to a class of obstacle problems with non-standard growth conditions.
The novelty here is that we are able to establish the local boundedness under a sharp bound   
on the gap between the growth exponents.
\end{abstract}


\section{Introduction}
This paper deals with the local boundedness of the solutions $u \in W^{1,p}(\Omega)$ to variational obstacle problems of the type
\begin{gather}\label{obpro}
\min \biggl\{ \io F(x,z,Dz) \dd x \ : \ z \in \mathcal{K}_{\psi}(\Omega)  \biggr\}.
\end{gather}
Here $\Omega$ is a bounded open set of $\mathbb{R}^n$, $n \geq 2$, the function $\psi: \Omega \rightarrow [-\infty, + \infty)$, called \textit{obstacle}, belongs to the Sobolev class $W^{1,p}(\Omega)$ and 
\begin{gather}\label{Insieme ammissibile}
\mathcal{K}_{\psi}(\Omega) := \{ z \in u_0 + W_0^{1,p}(\Omega) : z \geq \psi \ \text{a.e. in} \ \Omega  \} 
\end{gather}
is the class of admissible functions, where $u_0 \in W^{1,p}(\Omega)$ is a fixed boundary datum. To avoid trivialities, in what follows we shall assume that $\mathcal{K}_{\psi}(\Omega)$ is not empty.
We assume that the energy density $F: \Omega \times \R \times \R^n \rightarrow \R$ is a Carathéodory function such that 
\begin{gather}
(s,\xi) \mapsto F(x,s,\xi) \ \text{is convex},  \label{hp1} \\
c_1 |\xi|^p \leq F(x,s,\xi) \leq c_2 (1+ |s|^\gamma+|\xi|^q), \label{hp2}
\end{gather}
for almost all $x \in \Omega$ and all $s \in \R$, $\xi \in \R^n$, where $ 1 \leq p \leq q$, $0 \leq \gamma$ and $c_1,c_2 >0$ are fixed constants.  We recall that $u \in \mathcal{K}_{\psi}(\Omega)$ is a solution to \eqref{obpro} if and only if  $F(x,u, Du)\in L_{loc}^1(\Omega)$ and the minimality condition
\begin{equation} \label{Soluzione di minimo}
\int_{\supp(u-\varphi)}F(x,u,Du)\dd x \le \int_{\supp(u-\varphi)}F(x, \varphi,D\varphi)\dd x
\end{equation}
is satisfied for all $\varphi \in \mathcal{K}_{\psi}(\Omega)$.

The obstacle problem has been first considered in the works by Stampacchia \cite{Articolo16}, related to the capacity of a set in potential theory, and Fichera \cite{Articolo17}, who solved the so-called \textit{Signorini problem} in elastostatics.
\\
The study of the local boundedness for minimizers of integral functionals is a classic topic in Partial Differential Equations and Calculus of Variations starting from the
classical result by De Giorgi \cite{Libro1} and often is the first step in the analysis of the regularity properties of the solutions.
\\
In the last years there has been an intense research activity concerning the local boundedness of minimizers of unconstrained problems, in case of energy densities satisfying non-standard growth conditions (see for example \cite{Articolo14,Articolo6,Articolo15,Articolo9,Articolo13,Articolo11,Articolo12}). 
  The analogous study for solutions of obstacle problems under non standard growth has been exploited  in \cite{Articolo3}, where the local boundedness has been  established  assuming that 
  $$q <p_n^*=\frac{np}{n-p}$$
 provided the obstacle is locally bounded.
 We also point out a particular case of $p,q$-growth condition considered by Chlebicka and De Filippis \cite{Articolo18}, where the local boundedness for solutions to obstacle problems is proven for double phase energy densities.
 \\
  A restriction on the closeness between the growth exponents cannot be avoided, since  
    thanks to the well known Marcellini's counterexample \cite{Articolo4}, we know that if
\[
q>\frac{(n-1)p}{n-1-p}=p_{n-1}^* \hspace{1.2cm} 1<p<n-1
\]
minimizers of functionals with $p,q$-growth can be unbounded even in the unconstrained setting.
In a very recent paper  \cite{Articolo5},  Hirsch and Sch$\mathrm{\ddot{a}}$ffner proved that the sharp bound \eqref{gap}
    is sufficient to the local boundedness for unconstrained minimizers.  
  \\
    The aim of this paper is  to prove that \eqref{gap}
    is sufficient also to establish the local boundedness of solutions to the obstacle problem.

    More precisely, we are to going to prove the following

\begin{thm}\label{mainthm}
Let $u \in W^{1,p}(\Omega)$ be a solution to \eqref{obpro} under assumptions \eqref{hp1} and \eqref{hp2}, for exponents $1 \leq p < q $, $0 \leq \gamma$ verifying 
\begin{equation}\label{gap}
  	\frac{1}{q}\ge \frac{1}{p}-\frac{1}{n-1}
  \end{equation}
and
\begin{equation}\label{crescita in u}
    \quad \gamma \leq p^*_n:=
    \begin{cases}
\frac{np}{n-p} & \quad \text{if} \  p < n \\
\text{any finite exponent in } [q, \infty) & \quad \text{if} \ p \geq n
\end{cases}.
\end{equation}
If $\psi \in L^\infty_{\text{loc}}(\Omega)$, then $u \in L^\infty_{\text{loc}}(\Omega)$ and the following estimate
$$\sup_{B_{R_0/2}}|u| \leq C (\sup_{B_{R_0}}|\psi|+\Vert u \Vert_{W^{1,p}(B_{R_0})})^\pi,$$
holds for every ball $B_{R_0} \Subset \Omega$, for $\pi:= \pi(n,p,q)$ and with $C:=C(n,p,q,R_0)$.
\end{thm}
Observe that, due to the local nature of our regularity results, we are not requiring further properties on the boundary datum $u_0$ in \eqref{Insieme ammissibile}.\\ 
The proof of Theorem \ref{mainthm} is achieved following the strategy first proposed in \cite{Articolo3}, i.e. using the well known De Giorgi method that consists in deriving a suitable Caccioppoli inequality on the superlevel sets of the solution to \eqref{obpro}.
In order to do so, one has to use test functions obtained truncating the solution. Here, the difficulties come from the set of admissible test functions that must belong to the admissible class $\mathcal{K}_\psi(\Omega)$ and this is where the local boundedness of the obstacle $\psi$ comes into play.
We also remark that the crucial tool in order to achieve the result under the sharp bound on the gap between the exponents is the Sobolev inequality on the spheres as done in \cite{Articolo5}.
\\
 It is worth pointing out that assumption \eqref{gap} is essentially sharp in order to establish local boundedness for solutions of \eqref{obpro}. Indeed, in view of a counterexample by Franchi, Serapioni and Serra Cassano \cite{Articolo7}, the conclusion of Theorem \ref{mainthm} is false if condition \eqref{gap} is replaced by
\[
\frac{1}{p}+\frac{1}{q} < \frac{1}{n-1}+\varepsilon
\] 
for any $\varepsilon > 0$, already for unconstrained minimizers.
\\
The paper is organized as follows: in Section 2 we introduce some notations and collect some results that will be needed in the sequel; in Section 3 we derive a Caccioppoli inequality for the minimizer of \eqref{obpro}; Section 4 is devoted to the proof of Theorem \ref{mainthm}.

\section{Notations and Preliminary Results}
In this paper we will denote by $C$ or $c$ a general positive constant that may vary from line to line. Relevant dependencies on parameters will be highlighted using parentheses.
With the symbol $B_{r}(x)$ we will denote the ball with centre $x$ and radius $r$. We shall omit the dependence on the center and on the radius when they are clear from the context.
Moreover, we shall denote by $S_r = \{ x \in \R^n : |x| =r \}$ the sphere of radius $r$ on $\R^n$.

The following is a well-known iteration lemma (see \cite[Lemma 6.1]{Libro1} for the proof).
\begin{lem}\label{lm1}
Let $\Phi  :  [\rho,R] \rightarrow \mathbb{R}$ be a bounded nonnegative function. Assume that for all $\rho \leq t < s \leq R$ it holds
$$\Phi (t) \leq \theta \Phi(s) +A + \dfrac{B}{(s-t)^2}+ \dfrac{C}{(s-t)^{\gamma}}$$
where $\theta \in (0,1)$, $A$, $B$, $C \geq 0$ and $\gamma >0$ are constants. Then there exists a constant $c=c(\theta, \gamma)$ such that
$$\Phi (\rho ) \leq c \biggl( A+ \dfrac{B}{(R-\rho)^2}+ \dfrac{C}{(R-\rho)^{\gamma}}  \biggr).$$
\end{lem}

A key ingredient in the proof of Theorem \ref{mainthm} is the following lemma, that can be found in \cite[Lemma 2.1]{Articolo5}.
\begin{lem}\label{lm2}
Let $n \geq 2$. For any $0 < \rho < \sigma < \infty$, $v \in L^1(B_\sigma)$ and $s >1$, we set
$$
I(\rho, \sigma,v):= \inf \biggr\{ \int_{B_\sigma} |v||D\eta|^s \dd x : \eta \in \mathcal{C}^1_0(B_\sigma), \eta \geq 0, \eta=1 \ \text{in} \ B_\rho \biggl\}.
$$
Then for every $\delta \in (0,1]$
$$I(\rho,\sigma,v) \leq (\sigma-\rho)^{s-1+\frac{1}{\delta}}
\biggr( \int_\rho^\sigma \biggr( \int_{S_r}|v| \dd  \mathcal{H}^{n-1}  \biggl)^\delta \dd r\biggl)^{\frac{1}{\delta}}.$$
\end{lem}

Next Lemma, whose proof can be found in \cite[Lemma $7.1$]{Libro1}, allows us to iterate the Caccioppoli type estimate and it is crucial to establish the local boundedness result.

\begin{lem}\label{lm3}
Let $\alpha >0$ and let $(J_i)$ be a sequence of real positive numbers, such that
$$J_{i+1} \leq A \lambda^i J_i^{1+\alpha}, $$
with $A>0$ and $\lambda >1$. If $J_0 \leq A^{-\frac{1}{\alpha}} \lambda^{-\frac{1}{\alpha^2}}$, then $$J_i \leq \lambda^{-\frac{i}{\alpha}}J_0 \quad \text{and} \quad \lim_{i \rightarrow + \infty}J_i=0.$$
\end{lem}

We conclude this subsection recalling the Sobolev inequality on spheres (see e.g. \cite[Chapter 16]{Libro2}).

\begin{lem}\label{sobineq}
Let $v \in W^{1,m}(S_1, \dd \mathcal{H}^{n-1})$ with $m \in [1,n-1)$. Then there exists $c=c(n,m)$ such that
$$\biggr( \int_{S_1} |v|^{m^*} \dd \mathcal{H}^{n-1} \biggl)^\frac{1}{m^*} \leq c \biggr( \int_{S_1} (|D v|^m+|v|^m) \dd \mathcal{H}^{n-1} \biggl)^\frac{1}{m},$$
where $\frac{1}{m^*}=\frac{1}{m}-\frac{1}{n-1}$. 
\end{lem}

\section{Caccioppoli Inequality}
If $u \in W^{1,p}(\Omega)$, $k \in \R$ and $B_R \subset \Omega$ is a ball, we set
$$A_{k,R}:= \{ x \in B_R : u(x) > k \}.$$
The main result of this section is the following Caccioppoli inequality.

\begin{thm}\label{cacciopthm}
Let $u \in W^{1,p}(\Omega)$ be a solution to \eqref{obpro} under assumptions \eqref{hp1} and \eqref{hp2}, for exponents $1 \leq p < q $ verifying \eqref{gap} and $0 \leq \gamma$.
Assume that $\psi \in L^\infty_{\text{loc}}(\Omega)$. Then the following inequality 
\begin{align}\label{caccioppoli}
 \int_{B_\rho} |D(u-k)_+|^p \dd x 
\leq   
\dfrac{C}{(R-\rho)^\mu} \Vert (u-k)_+\Vert^q_{W^{1,p}(B_R)} |A_{k,R}|^{q(\frac{1}{q_*}-\frac{1}{p})} + C k^\gamma |A_{k,R}|
\end{align}
holds for all balls $B_\rho \subset B_R \subset B_{R_0} \Subset \Omega$ and for every $k \geq \max \{ \sup_{B_{R_0}} \psi , 1 \}$, where  
\begin{equation}\label{sobq}
    \frac{1}{q_*}:= \min \biggl\{ \frac{1}{q}+\frac{1}{n-1},1  \biggr\} ,
\end{equation}
$\mu := q-1+\frac{q}{q_*}$ and with $C:=C(n,q,R_0)$.
\end{thm}
\proof
Let us fix $B_{R_0} \Subset \Omega$. Let $\frac{R_0}{2} \leq \rho \leq s < t \leq R \leq R_0 $ and let $\eta \in \mathcal{C}_0^\infty(B_t)$ be a cut-off function such that $0 \leq \eta \leq 1$, $\eta=1 $ in $B_s$, $|D \eta | \leq \frac{2}{t-s}$. By virtue of the assumption $\psi \in L^\infty_{\text{loc}}(\Omega)$, we may fix $k \geq \max \{ \sup_{B_{R_0}} \psi , 1 \}$ and define $u_k := (u-k)_+=\max \{ u-k,0 \}$. Note that $\varphi = u- \eta^{\sigma}u_k$, with $\sigma >q$, belongs to   $\K$. Indeed,
\begin{equation}
\varphi =
\begin{cases}
u \geq \psi & \quad \text{if} \ u \leq k \\
u- \eta^{\sigma}(u-k)=(1-\eta^\sigma)(u-k)+k \geq k \geq \psi & \quad \text{if} \ u \geq k
\end{cases}.
\end{equation}
By the minimality of $u$ and the convexity assumption \eqref{hp1}, we get
\begin{align*}
& \int_{A_{k,s}} F(x,u,Du) \dd x \leq \int_{A_{k,t}} F(x,u+ \varphi, Du+ D \varphi) \dd x\\
& = \int_{A_{k,t}}F\biggl(x,(1-\eta^\sigma)u+\eta^\sigma k, (1-\eta^\sigma)Du-\sigma \eta^\sigma \dfrac{D \eta}{\eta}(u-k) \biggr) \dd x \\
& \leq \int_{A_{k,t}} \biggl[(1-\eta^\sigma)F(x,u,Du) +\eta^\sigma F \biggl( x,k, -\sigma \frac{D \eta}{\eta}(u-k) \biggr)\biggr]\dd x.
\end{align*}
Taking into account that $\text{supp} (1- \eta^\sigma) \subset A_{k,t} \setminus A_{k,s}$, $\sigma >q$, $t \leq R$ and the growth assumption \eqref{hp2}, we obtain
\begin{align}
& \int_{A_{k,s}} F(x,u,Du) \dd x \notag\\
& \leq \int_{A_{k,t} \setminus A_{k,s}} (1-\eta^\sigma)F(x,u,Du) \dd x + c_2
\int_{A_{k,t}} \eta^\sigma \biggl( 1+k^q+\sigma^q\frac{|D \eta |^q}{\eta^q}(u-k)^q \biggr) \dd x  \notag\\
& \leq \int_{A_{k,t} \setminus A_{k,s}} F(x,u,Du) \dd x +
C k^q |A_{k,R}|
+C \int_{A_{k,t}}|D \eta|^q (u-k)^q \dd x . \label{3.0}
\end{align}
Exploiting Lemma \ref{sobineq} after the use of H\"{o}lder's inequality in case $q \in (1, \frac{n-1}{n-2})$, we infer
\begin{equation}\label{sobsphere}
\biggr( \int_{S_1} |v|^{q} \dd \mathcal{H}^{n-1} \biggl)^\frac{1}{q} \leq c(n,q) \biggr( \int_{S_1} (|D v|^{q_*}+|v|^{q_*}) \dd \mathcal{H}^{n-1} \biggl)^\frac{1}{q_*},
\end{equation}
for every $v \in W^{1,q_*}(S_1, \dd \mathcal{H}^{n-1})$, where $q_* \geq 1$ is given by 
$$
\frac{1}{q_*}:= \min \biggr\{ \frac{1}{q}+\frac{1}{n-1},1  \biggl\}. $$
\\We set
$$\mathcal{A}(s,t) := \{ \eta \in \mathcal{C}^\infty_0(B_t): \eta =1 \ \text{in} \ B_s  \}.$$
Combining \eqref{sobsphere} applied to $v(y)=u^q_k(ry)$ with $r \in (s,t)$ and Lemma \ref{lm2} with $\delta= \frac{q_*}{q}$, we get
\begin{align*}
& \inf_{\mathcal{A}(s,t)} \int_{B_t} |D \eta|^q u_k^q \dd x \\
& \leq (t-s)^{-(q-1+\frac{q}{q_*})} \biggr( \int_s^t r^{(n-1){\frac{q_*}{q}}}\biggr( \int_{S_1}  u^q_k(ry) \dd \mathcal{H}^{n-1}(y) \biggl)^{\frac{q_*}{q}}\dd r \biggl)^{\frac{q}{q_*}}\\
& \leq c(n,q)(t-s)^{-(q-1+\frac{q}{q_*})} \biggr( \int_s^t r^{(n-1){\frac{q_*}{q}}}\int_{S_1} (r^{q_*} |Du_k(ry)|^{q_*} \\& \quad +u_k^{q_*}(ry))\dd \mathcal{H}^{n-1}(y) \dd r \biggl)^{\frac{q}{q_*}} \\
& = c(n,q)(t-s)^{-(q-1+\frac{q}{q_*})} \biggr( \int_s^t r^{(n-1)q_*(\frac{1}{q}+\frac{1}{n-1}-\frac{q_*}{q})}\int_{S_r} ( |Du_k(y)|^{q_*} \\
& \quad +r^{-q_*}u_k^{q_*}(y))\dd \mathcal{H}^{n-1}(y) \dd r \biggl)^{\frac{q}{q_*}}.
\end{align*}
Using the fact that $\frac{R_0}{2} \leq s < t \leq R_0$, we deduce 
\begin{align}
& \inf_{\mathcal{A}(s,t)} \int_{B_t} |D \eta|^q u_k^q \dd x \notag \\
& \leq C(n,q,R_0)(t-s)^{-(q-1+\frac{q}{q_*})} \biggr( \int_s^t \int_{S_r} ( |Du_k(y)|^{q_*} +u_k^{q_*}(y))\dd \mathcal{H}^{n-1}(y) \dd r \biggl)^{\frac{q}{q_*}} \label{3.1}.
\end{align}
Now, notice that inequality \eqref{gap} implies $q_* \leq p$. Thus, using H\"{o}lder's inequality in estimate \eqref{3.1}, we obtain
\begin{align}
 \inf_{\mathcal{A}(s,t)} \int_{B_t} |D \eta|^q u_k^q \dd x \leq & \dfrac{C(n,q,R_0)}{(t-s)^\mu} \Vert u_k\Vert^q_{W^{1,q_*}(B_t \setminus B_s)} \notag \\
\leq & \dfrac{C(n,q,R_0)}{(t-s)^\mu} \Vert u_k\Vert^q_{W^{1,p}(B_t \setminus B_s)} |A_{k,t}|^{q(\frac{1}{q_*}-\frac{1}{p})}\notag\\
 \leq & \dfrac{C(n,q,R_0)}{(t-s)^\mu} \Vert u_k\Vert^q_{W^{1,p}(B_R)} |A_{k,R}|^{q(\frac{1}{q_*}-\frac{1}{p})},
 \label{3.2}
\end{align}
where we denote $\mu = q-1+\frac{q}{q_*}$.\\
Since estimate \eqref{3.0} holds for every $\eta \in \mathcal{A}(s,t)$, by \eqref{3.2}, we get
\begin{align*}
& \int_{A_{k,s}} F(x,u,Du) \dd x \notag\\
& \leq \int_{A_{k,t} \setminus A_{k,s}} F(x,u,Du) \dd x +
 C k^\gamma|A_{k,R}|  \notag\\
& \quad +\dfrac{C(n,q,R_0)}{(t-s)^\mu} \Vert u_k\Vert^q_{W^{1,p}(B_R)} |A_{k,R}|^{q(\frac{1}{q_*}-\frac{1}{p})}.
\end{align*}
Adding the integral $\int_{A_{k,s}} F(x,u,Du) \dd x$ to both sides of the previous estimate, by Lemma \ref{lm1} we get
\begin{align}\label{3.3}
& \int_{A_{k,\rho}} F(x,u,Du) \dd x \notag\\
&\leq C(n,q,R_0) \biggr\{k^\gamma |A_{k,R}|  
+\dfrac{1}{(R-\rho)^\mu} \Vert u_k\Vert^q_{W^{1,p}(B_R)} |A_{k,R}|^{q(\frac{1}{q_*}-\frac{1}{p})}  \biggl\}.
\end{align}
Eventually, inequality \eqref{3.3} and the growth condition at \eqref{hp2} yield
\begin{align*}
& \int_{B_\rho} |D(u-k)_+|^p \dd x  \notag \\
&\leq  C(n,q,R_0) \biggr\{ k^\gamma |A_{k,R}|  
+\dfrac{1}{(R-\rho)^\mu} \Vert u_k\Vert^q_{W^{1,p}(B_R)} |A_{k,R}|^{q(\frac{1}{q_*}-\frac{1}{p})}   \biggl\} ,
\end{align*}
i.e. the conclusion.
\endproof

\section{Proof of the Main Result}
Before giving the proof of Theorem \ref{mainthm}, we need to introduce some notations.\\
For a fixed ball $B_{R_0} \Subset \Omega$, define two sequences by setting
$$\rho_i :=\dfrac{R_0}{2} \biggr( 1+\dfrac{1}{2^i} \biggl),$$
$$k_i:= 2d \biggr( 1-\dfrac{1}{2^{i+1}}  \biggl),$$
where $d \geq \max \{ \sup_{B_{R_0}} \psi, 1 \}$ will be determined later. Note that 
$$\dfrac{R_0}{2} \leq \rho_i \leq R_0, \quad \rho_i \searrow \dfrac{R_0}{2}$$
and
$$d \leq k_i \leq 2d, \quad k_i \nearrow 2d.$$ 
The use of estimate \eqref{caccioppoli} on the concentric balls $B_{\rho_{i+1}} \subset B_{\rho_i}$ translates into
\begin{align}
&\Vert D(u-k_{i+1})_+ \Vert^p_{L^p(B_{\rho_{i+1}})} \notag\\
&\leq   
\dfrac{C}{(\rho_{i}-\rho_{i+1})^\mu} \Vert (u-k_{i+1})_+\Vert^q_{W^{1,p}(B_{\rho_{i}})} |A_{k_{i+1},\rho_i}|^{q(\frac{1}{q_*}-\frac{1}{p})} 
+ C k_{i+1}^\gamma|A_{k_{i+1},\rho_i}| , \label{caccioppoli1}
\end{align}
for every $i \in \N$.\\
Moreover, define the sequence $(J_i)$ setting
$$J_i:=\Vert (u-k_i)_+ \Vert^p_{W^{1,p}(B_{\rho_i})}.$$
\\We begin proving an inequality that will be crucial for the proof of Theorem \ref{mainthm}.
\begin{lem}\label{dis}
Let $u \in W^{1,p}(\Omega)$ be a solution to \eqref{obpro} under assumptions \eqref{hp1} and \eqref{hp2}, for exponents $1 \leq p < q $, $0 \leq \gamma$ verifying \eqref{gap}, \eqref{crescita in u} respectively. Then, there exists a constant $\tilde{C}>0$ such that for every $i \in \N$
\begin{equation}
J_{i+1} \leq \dfrac{\tilde{C}}{d^\theta} \lambda^{i}  J_i^{1+\alpha}, \notag
\end{equation}
for some $\lambda >1$ and $\alpha >0$, where $\dfrac{1}{d^\theta}:= \max \biggl\{ \dfrac{1}{d^{\frac{p^*_n p}{n}}},
\dfrac{1}{d^{q  (\frac{1}{q_*}-\frac{1}{p})p^*_n}}, \dfrac{1}{d^{p^*_n -\gamma}} \biggr\}$.
\end{lem}
\proof
We observe that $k_{i+1}-k_{i} < u-k_i$ on $A_{k_{i+1},\rho_{i}}$ for every $i \in \N$. By Sobolev inequality we get
\begin{align}
|A_{k_{i+1},\rho_{i}}| &  \leq\int_{A_{k_{i+1},\rho_{i}}} \biggr(\dfrac{u-k_i}{k_{i+1}-k_i}\biggl)^{p^*_n} \dd x 
 \leq \dfrac{\Vert (u-k_i)_+ \Vert^{p^*_n}_{L^{p^*_n}(B_{\rho_i})}}{(k_{i+1}-k_i)^{p^*_n}} \notag\\
& \leq C(n,p)\dfrac{J_i^{\frac{p^*_n}{p}}}{(k_{i+1}-k_i)^{p^*_n}}  , \label{est1}
\end{align}
where $p^*_n$ was introduced in \eqref{crescita in u}.\\
Furthermore, we have
\begin{align}
&  \Vert (u-k_{i+1})_+  \Vert_{W^{1,p}(B_{\rho_{i}})}^p \notag\\
& = \int_{B_{\rho_{i}}}(u-k_{i+1})_+^p \dd x  + \int_{B_{\rho_{i}}}|D(u-k_{i+1})_+|^p \dd x  \notag\\
&= \int_{A_{k_{i+1},\rho_{i}}}(u-k_{i+1})^p \dd x  + \int_{A_{k_{i+1},\rho_{i}}}|D(u-k_{i+1})|^p \dd x  \notag\\
& \leq \int_{A_{k_{i},\rho_{i}}}(u-k_{i})^p \dd x  + \int_{A_{k_{i},\rho_{i}}}|D(u-k_{i})|^p \dd x = J_i. \label{est2}
\end{align}
Inserting estimates \eqref{est1} and \eqref{est2} in \eqref{caccioppoli1}, we obtain
\begin{align}
\Vert D(u-k_{i+1})_+ \Vert^p_{L^p(B_{\rho_{i+1}})} \leq   
\dfrac{C}{(\rho_{i}-\rho_{i+1})^\mu} \dfrac{J_i^{\frac{q}{p}\bigl(1+p^*_n(\frac{1}{q_*}-\frac{1}{p}) \bigr)}}{(k_{i+1}-k_i)^{q  (\frac{1}{q_*}-\frac{1}{p})p^*_n}} + C \dfrac{k_{i+1}^\gamma J_i^{\frac{p^*_n}{p}}}{(k_{i+1}-k_i)^{p^*_n}}. \label{est3}
\end{align}
By using H\"{o}lder's and Sobolev inequalities and estimate \eqref{est1}, it follows
\begin{align}
\Vert (u-k_{i+1})_+  \Vert^p_{L^p(B_{\rho_{i+1}})}
\leq \Vert (u-k_{i})_+  \Vert^p_{L^{p^*_n}(B_{\rho_{i}})}
|A_{k_{i+1},\rho_i}|^{\frac{p}{n}}
\leq C\dfrac{J_i^{\frac{p^*_n}{n}+1}}{(k_{i+1}-k_i)^{\frac{p^*_n p}{n}}}. \label{est4}
\end{align}
Summing \eqref{est3} and \eqref{est4}, we obtain
\begin{align}
J_{i+1} \leq C \biggl[ \dfrac{J_i^{\frac{p^*_n}{n}+1}}{(k_{i+1}-k_i)^{\frac{p^*_n p}{n}}}
+  \dfrac{1}{(\rho_{i}-\rho_{i+1})^\mu} \dfrac{J_i^{\frac{q}{p} \bigl(1+p^*_n(\frac{1}{q_*}-\frac{1}{p}) \bigr)}}{(k_{i+1}-k_i)^{q  (\frac{1}{q_*}-\frac{1}{p})p^*_n}} +  \dfrac{k_{i+1}^\gamma J_i^{\frac{p^*_n}{p}}}{(k_{i+1}-k_i)^{p^*_n}}
\biggr], \label{est5}
\end{align}
for a constant $C:=C(n,p,q,R_0)$.
\\Recalling the definition of $k_i$, $\rho_i$ and using the fact that $k_i \leq 2d$, we can write inequality \eqref{est5} as follows
\begin{align}
&J_{i+1} \notag\\
& \leq C \biggl[ \dfrac{J_i^{\frac{p^*_n}{n}+1}}{d^{\frac{p^*_n p}{n}}}2^{(i+2)\frac{p^*_n p}{n}}
+  \dfrac{2^{(i+2) \bigl(\mu+ q  (\frac{1}{q_*}-\frac{1}{p})p^*_n \bigr)}}{R^\mu_0} \dfrac{J_i^{\frac{q}{p} \bigl(1+p^*_n(\frac{1}{q_*}-\frac{1}{p}) \bigr)}}{d^{q  (\frac{1}{q_*}-\frac{1}{p})p^*_n}} +  \dfrac{2^{(i+2)p^*_n}}{d^{p^*_n - \gamma}}
J_i^{\frac{p^*_n}{p}}
\biggr] \notag\\
& \leq \dfrac{C}{R^\mu_0 d^\theta} 2^{i \tau}
(J_i^{\frac{p^*_n}{n}+1} +
J_i^{\frac{q}{p} \bigl(1+p^*_n(\frac{1}{q_*}-\frac{1}{p})\bigr)} + J_i^{\frac{p^*_n}{p}} ), \label{est6}
\end{align}
where
$$\dfrac{1}{d^\theta}:= \max \biggl\{ \dfrac{1}{d^{\frac{p^*_n p}{n}}},
\dfrac{1}{d^{q  (\frac{1}{q_*}-\frac{1}{p})p^*_n}}, \dfrac{1}{d^{p^*_n -\gamma}}
 \biggr\},$$
$$ \tau:= \max \biggl\{ \frac{p^*_n p}{n},  \mu+ q  \biggl(\frac{1}{q_*}-\frac{1}{p} \biggr)p^*_n, p^*_n \biggr\} .$$
Notice that by the definition of $J_i$ and $d$, we easily derive 
$$J_i \leq \Vert (u-\sup_{B_{R_0}} \psi)_+ \Vert^p_{W^{1,p}(B_{R_0})} \quad \forall i \in \N.$$
If $q \leq p^*_n $, setting
$$\beta := \max \biggr\{ \frac{p^*_n}{n}+1 - \dfrac{q}{p},
 \frac{q }{p} \biggl(\frac{1}{q_*}-\frac{1}{p} \biggr) p^*_n,
 \frac{p^*_n}{p} - \dfrac{q}{p} \biggl\} > 0,$$
estimate \eqref{est6} can be written as
\begin{align*}
J_{i+1} 
& \leq \dfrac{C}{R^\mu_0 d^\theta} 2^{i \tau}
(1+ \Vert (u-\sup_{B_{R_0}} \psi)_+ \Vert^p_{W^{1,p}(B_{R_0})})^\beta J_i^{\frac{q}{p}} 
=\dfrac{\tilde{C}}{d^\theta} 2^{i \tau}  J_i^{\frac{q}{p}} .
\end{align*}
Now, suppose $q > p^*_n$.
In this case $p<n$ (see \eqref{crescita in u} for the definition of $p^*_n$). Thus, $\frac{p^*_n}{n}+1= \frac{p^*_n}{p}$.
Setting
$$\sigma := 
 \frac{q}{p} \biggl(1+ \biggl(\frac{1}{q_*}-\frac{1}{p} \biggr) p^*_n \biggr) - \frac{p^*_n}{p} >0,$$
estimate \eqref{est6} can be written as
\begin{align*}
J_{i+1} 
& \leq \dfrac{C}{R^\mu_0 d^\theta} 2^{i \tau}
(1+ \Vert (u-\sup_{B_{R_0}} \psi)_+ \Vert^p_{W^{1,p}(B_{R_0})})^\sigma J_i^{\frac{p^*_n}{p}} 
=\dfrac{\tilde{C}}{d^\theta} 2^{i \tau}  J_i^{\frac{p^*_n}{p}} .
\end{align*}
\endproof

Now, we are in position to give the proof of our main result.\\
{\it Proof of Theorem \ref{mainthm}}
Notice that $u \geq \psi$ and $\psi \in L^\infty_{\text{loc}}(\Omega)$ imply that $u$ is bounded from below. Hence, we only need to show that $u$ is bounded from above.\\
In order to show that the minimizer $u$ is bounded in the ball $B_{R_0/2}$, we prove that the sequence $(J_i)$ satisfies the assumptions of Lemma \ref{lm3}. 
Indeed, Lemma \ref{lm3} yields 
\begin{equation}
\lim_{i \rightarrow + \infty}J_i=0.  \label{lim0}
\end{equation}
On the other hand, we have
\begin{equation}\label{lim1}
\lim_{i \rightarrow + \infty}J_i = \Vert (u-2d)_+ \Vert^p_{W^{1,p}(B_{R_0/2})}. 
\end{equation}
Thus, \eqref{lim0} and \eqref{lim1} allow us to conclude that
$$\sup_{B_{R_0/2}}u \leq 2d.$$
Denote
$$\varepsilon := \frac{1}{q_*}-\frac{1}{p} .$$
If $\gamma < p^*_n$ and $\varepsilon >0$, by Lemma \ref{dis}
the sequence $(J_i)$ satisfies assumptions of Lemma \ref{lm3} with
$$A= \dfrac{\tilde{C}}{d^\theta} , \quad
\lambda= 2^\tau \quad \text{and} \quad \alpha= \dfrac{q}{p}-1 \ \text{or} \   \alpha= \dfrac{p^*_n}{p}-1.$$
We also have that 
$$J_0 \leq A^{-\frac{1}{\alpha}}\lambda^{-\frac{1}{\alpha^2}}\quad \text{with} \quad d \geq \tilde{C}^{\frac{1}{\theta}} \lambda^{\frac{1}{\alpha \theta}} \Vert u \Vert^{p \frac{\alpha}{\theta}}_{W^{1,p}(B_{R_0})} .$$
Therefore, \eqref{lim0} holds.\\
If $\gamma<p^*_n$ and $\varepsilon =0$, 
Lemma \ref{dis} implies that the sequence $(J_i)$ satisfies Lemma \ref{lm3} with
$$A= \tilde{C} , \quad
\lambda= 2^\tau \quad \text{and} \quad \alpha= \dfrac{q}{p}-1 \ \text{or} \   \alpha= \dfrac{p^*_n}{p}-1 .$$
Therefore, by Lemma \ref{lm3} we have that \eqref{lim0} holds if
\begin{equation}
J_0 \leq A^{-\frac{1}{\alpha}} \lambda^{-\frac{1}{\alpha^2}} .  \label{jo}
\end{equation}
By definition, $J_0 = \Vert (u-d)_+ \Vert^p_{W^{1,p}(B_{R_0})}$.
We choose $d>0$ large enough so that \eqref{jo} holds; this is possible since $u \in W^{1,p}(B_{R_0})$ and 
$$J_0 = \Vert (u-d)_+ \Vert^p_{W^{1,p}(A_{d,R_0})} \leq 
\Vert u \Vert^p_{W^{1,p}(A_{d,R_0})}  \longrightarrow 0 \quad \text{as} \ d \rightarrow + \infty.$$
With this choice of $d$, we get \eqref{lim0}.\\
Eventually, the case $\gamma=p^*_n$ and $\varepsilon \geq 0$ can be treated as the previous one.
\qed
\vspace{0.5cm}
\\
\textbf{Remark 4.2.}
Let us consider functionals $F$ satisfying the double-side bound
\begin{align}\label{crescitadf}
    \nu H(x,z) \leq F(x,s,z) \leq L H(x,z),
\end{align}
for constants $0 \leq \nu \leq L$, where
\begin{align}
    H(x,z):= |z|^p+a(x)|z|^q,
\end{align}
with $1 \leq p < q$ and $0 \leq a(x) \in L^\infty(\Omega)$. Then, under the closeness condition \eqref{gap} on the exponents $p,q$, solutions to \eqref{obpro} with $F$ verifying \eqref{crescitadf} are locally bounded provided the obstacle is locally bounded.
The arguments are essentially the same of Theorem \ref{mainthm}. Indeed, an analogous version of Theorem \ref{cacciopthm} and Lemma \ref{dis} can be proved similarly with exponent $\gamma =0$.
\\Let $u \in W^{1,p}(\Omega)$ be a solution to \eqref{obpro} with energy density satisfying the growth condition \eqref{crescitadf}.
Using the same notation introduced in the proof of Theorem \ref{cacciopthm}, we have
\begin{align} 
& \int_{A_{k,s}} F(x,u,Du) \dd x \leq \int_{A_{k,t}} F(x,u+ \varphi, Du+ D \varphi) \dd x \notag\\
& \leq L \int_{A_{k,t}} H (x, Du+D \varphi ) \dd x  \notag\\
& = L\int_{A_{k,t}} H \biggl(x, (1-\eta^\sigma)Du- \sigma \eta^\sigma \frac{D \eta}{\eta} (u-k) \biggr) \dd x. \label{6.1}
\end{align}
By the very definition of $H(x,z)$, we have that for a.e. $x \in \Omega $ and every $z_1,z_2 \in \R^n$
\begin{equation}
    H(x,z_1+z_2) \leq 2^q (H(x,z_1)+H(x,z_2)  ). \label{6.0}
\end{equation}
Using \eqref{6.0}, the fact that $\eta \leq 1$ and the boundedness of function $a$, we estimate inequality \eqref{6.1} as follows
\begin{align*} 
\int_{A_{k,s}} & F(x,u,Du) \dd x \\
\leq  & C(q,L)\int_{A_{k,t}} H (x, (1-\eta^\sigma)Du ) \dd x \\
&+  C(q,L)\int_{A_{k,t}} H \biggl(x, - \sigma \eta^\sigma \frac{D \eta}{\eta} (u-k) \biggr) \dd x \\
\leq &  C(q,L)\int_{A_{k,t}} (1-\eta^\sigma)^p H (x, Du ) \dd x \\
&+  C(q,L, \Vert a \Vert_\infty)\int_{A_{k,t}}  (1+|D \eta|^q (u-k)^q ) \dd x .
\end{align*}
Taking into account that $\text{supp} (1- \eta^\sigma) \subset A_{k,t} \setminus A_{k,s}$, the growth assumption \eqref{crescitadf} and $t \leq R$, we obtain
\begin{align*}
& \int_{A_{k,s}} H(x,Du) \dd x 
  \leq C \int_{A_{k,s}} F(x,u,Du) \dd x \notag\\
& \leq C \int_{A_{k,t} \setminus A_{k,s}} H(x,Du) \dd x+
 C |A_{k,R}|  
+C \int_{A_{k,t}}|D \eta|^q (u-k)^q \dd x . 
\end{align*}
Recalling the definition of $H(x,Du)$, the previous estimate yields the Caccioppoli inequality
\eqref{caccioppoli} with $\gamma=0$.
\\It is worth noting that assumption \eqref{gap} improves the bound on the gap $q/p$ established in \cite{Articolo18} for obtaining the local boundedness of solutions to double phase obstacle problems.

\section{Conclusions}
The main result of this paper states that solutions to a class of obstacle problems, satisfying $p,q$-growth conditions under the sharp relation
\[	
\frac{1}{q}\ge \frac{1}{p}-\frac{1}{n-1},
\]
are locally bounded.\\
Let us briefly describe the contents of the paper. After recalling some notation and preliminary results in Section 2, we focus on deriving the intermediate steps that will put us in the position to prove our main result. In particular, in Section 3, we show that the solution to the obstacle problem satisfies a suitable Caccioppoli inequality (see Theorem \ref{cacciopthm}). In Section 4, we give the proof of Lemma \ref{dis} that allows us to iterate the Caccioppoli type estimate, so that we are eventually able to prove Theorem \ref{mainthm}.

\end{document}